\documentclass{amsart}
\usepackage{amssymb,latexsym,amsmath,amsfonts,hhline}
\usepackage{pdfsync,comment}
\usepackage[usenames]{color} 

\input xy
\xyoption{all}

\newcommand{\cal}{\mathcal}

\newtheorem{formula}{}[section]
\newtheorem{definition}[formula]{Definition}
\newtheorem{corollary}[formula]{Corollary}
\newtheorem{proposition}[formula]{Proposition}
\newtheorem{remark}[formula]{Remark}
\newtheorem{lemma}[formula]{Lemma}
\newtheorem{theorem}[formula]{Theorem}

\def\thrm{\begin{theorem}}
\def\thrml#1{\begin{theorem}\label{#1}}
\def\ethrm{\end{theorem}}
\def\rmrk{\begin{remark}}
\def\rmrkl#1{\begin{remark}\label{#1}}
\def\ermrk{\end{remark}}
\def\dfntn{\begin{definition}}
\def\dfntnl#1{\begin{definition}\label{#1}}
\def\edfntn{\end{definition}}
\def\nmrt{\begin{enumerate}}
\def\enmrt{\end{enumerate}}
\def\tm#1{\item[{\rm (#1)}]}
\def\qtn{\begin{equation}}
\def\qtnl#1{\begin{equation}\label{#1}}
\def\eqtn{\end{equation}}
\def\lmm{\begin{lemma}}
\def\lmml#1{\begin{lemma}\label{#1}}
\def\elmm{\end{lemma}}
\def\crllr{\begin{corollary}}
\def\crllrl#1{\begin{corollary}\label{#1}}
\def\ecrllr{\end{corollary}}
\def\prpstn{\begin{proposition}}
\def\prpstnl#1{\begin{proposition}\label{#1}}
\def\eprpstn{\end{proposition}}
\def\css{\begin{cases}}
\def\ecss{\end{cases}}

\def\proof{\noindent{\bf Proof}.\ }

\def\cA{{\cal A}}

\def\cS{{\cal S}}

\def\fG{{\mathfrak G}}

\def\mF{{\mathbb F}}
\def\mZ{{\mathbb Z}}

\DeclareMathOperator{\aut}{Aut}
\DeclareMathOperator{\alt}{Alt}

\DeclareMathOperator{\cyc}{Cyc}

\DeclareMathOperator{\GF}{GF}
\DeclareMathOperator{\GL}{GL}

\DeclareMathOperator{\im}{im}

\DeclareMathOperator{\iso}{Iso}

\DeclareMathOperator{\orb}{Orb}

\DeclareMathOperator{\rad}{rad}

\DeclareMathOperator{\rk}{rk}

\DeclareMathOperator{\Span}{Span}

\DeclareMathOperator{\sym}{Sym}

\def\eprf{\hfill$\square$}

\def\qaq{\quad\text{and}\quad}
\def\qoq{\quad\text{or}\quad}

\newcommand{\grp}[1]{\langle {#1}\rangle}
\newcommand{\und}[1]{{\underline{#1}}}
\def\wh{\widehat}

\begin{document}
\title{Abelian Schur groups of odd order}

\author{Ilia Ponomarenko}
\address{St.Petersburg Department of the Steklov Mathematical Institute, St.Petersburg, Russia}
\thanks{The work of the first author was supported by the RFBR Grant No. 17-51-53007 GFEN\_a and by the RAS Program of Fundamental Research ``Modern Problems of Theoretical Mathematics''. 
The second author was  supported by RSF (project No. 14-21-00065)}
\email{inp@pdmi.ras.ru}
\author{Grigory Ryabov}
\address{Sobolev Institute of Mathematics, Novosibirsk, Russia}
\address{Novosibirsk State University, Novosibirsk, Russia}
\email{gric2ryabov@gmail.com}
\date{}

\begin{abstract}
A finite group $G$ is called a Schur group if any Schur ring over~$G$ is
associated in a natural way with a subgroup of $\sym(G)$ that contains all right translations.	It is proved that the group $C_3\times C_3\times C_p$ is Schur for any prime~$p$. Together with earlier results, this completes 
a classification of the abelian Schur groups of odd order.
\end{abstract}

\maketitle
\section{Introduction}

A {\it Schur ring} or {\it S-ring} over a finite group $G$ can be defined as a subring of the group ring $\mZ G$ that is a free $\mZ$-module spanned by a partition of $G$ closed under taking inverse and containing the identity element $e$ of $G$ as a class (see Section~\ref{031017a} for details). An important example of such a partition is given by the orbits of the point stabilizer $K_e$ of a permutation group $K$ such that
\qtnl{031017u}
G_{right}\le K\le\sym(G),
\eqtn
where $G_{right}$ is the group induced by the right translations of~$G$. The corresponding S-rings are said to be {\it schurian} in honor of I.~Schur who studied the S-rings of this type.\medskip 

In fact, there are a lot of non-shurian S-rings. An infinite family of them can be found in paper of R.~P\"oschel \cite{Po}, where he introduced a concept of a {\it Schur group}: the group $G$ is Schur if every S-ring over $G$ is   schurian. A motivation for being interested in the Schur  groups comes from the problem of testing isomorphism of Cayley graphs, see~\cite{KP81,NP}.\medskip

In~\cite{EKP2}, it was proved that every finite abelian Schur group belongs to one of several explicitly given families. Recent results \cite{MP3,R17} show that two of them are indeed consist of Schur groups. The main result of the present paper is given in the theorem below concerning the next family.

\thrml{230917a}
For any prime $p$, all S-rings over a group $G=E_9\times C_p$ are schurian. In particular, $G$ is a Schur group.
\ethrm

All cyclic Schur groups were classified in \cite{EKP1}. Therefore, as an immediate consequence of this theorem and the above mentioned results, we obtain a classification of all abelian Schur groups of odd order.

\thrm
A noncyclic abelian  group of odd order is Schur if and only if it is  isomorphic to $C_3\times C_{3^k}$ for an integer $k\ge 1$, or $E_9\times C_p$ for a prime $p\ge 3$.
\ethrm

In Section~\ref{031017i}, we deduce Theorem~\ref{230917a} from Theorem~\ref{main} stating that any S-ring~$\cA$ over the group $G=E_9\times C_p$ is either obtained from two S-rings over smaller groups (and then the schurity of $\cA$ is under control) or is a cyclotomic S-ring (and then $\cA$ is schurian by definition). The proof of Theorem~\ref{main} is mainly based on Theorem~\ref{240917a} giving a sufficient condition for $\cA$ to be cyclotomic. This carried out in three steps. First, the S-rings over $E_9$ are completely described (Section~\ref{020917f}). This enables us to prove that any class of the partition of~$G$ associated with~$\cA$ is an orbit of a suitable subgroup of~$\aut(G)$ (Section~\ref{130917a}). At the last step, we show that this subgroup can be chosen the same for all classes (Section~\ref{270917a}). For reader's convenience, we cite basic facts on S-rings in Section~\ref{031017a}.\medskip

{\bf Notation.}

As usual by $\mZ$ we denote the ring of rational integers.

A finite field of order $q$ is denoted by $\GF_q$.

The projections of $X\subseteq A\times B$ to $A$ and $B$ are denoted by $X_A$ and $X_B$, respectively.

The set of non-identity elements of a group $G$ is denoted by  $G^\#$.

The center of a group $G$ is denoted by $Z(G)$.

Let $X\subseteq G$. The subgroup of $G$ generated by $X$ is denoted by $\grp{X}$; we also set $\rad(X)=\{g\in G:\ gX=Xg=X\}$.

Let $\sigma\in\aut(G)$. The element $\sum_{x\in X}x^\sigma$ of the group ring $\mZ G$ is denoted by $\und{X}^\sigma$, and by $\und{X}$ if $\sigma$ is the identity.

The componentwise multiplication in the ring $\mZ G$ is denoted by $\circ$.

The group of all permutations of $G$ is denoted by $\sym(G)$.

The induced action of $G\le\sym(\Omega)$ on an invariant set $\Delta\subseteq\Omega$ is denoted by $G^\Delta$.

The cyclic group of order $n$ is denoted by  $C_n$.

The elementary abelian group of order $p^k$ is denoted by $E_{p^k}$.

\section{Preliminaries}\label{031017a}

\subsection{Definitions}
Let $G$ be a finite group. A subring~$\cA$ of the group ring~$\mZ G$ is  called a {\it Schur ring} ({\it S-ring}, for short) over~$G$ if there exists a partition $\cS=\cS(\cA)$ of~$G$ such that
\nmrt
\tm{S1} $\{1_G\}\in\cS$,
\tm{S2} $X\in\cS\ \Rightarrow\ X^{-1}\in\cS$,
\tm{S3} $\cA=\Span\{\und{X}:\ X\in\cS\}$.
\enmrt
In particular, condition (S3) implies that $\cA$ is closed with respect to the componentwise multiplication~$\circ$. A group isomorphism $f:G\to G'$ is called a {\it Cayley isomorphism} from an S-ring $\cA$ over $G$ onto an S-ring $\cA'$ over $G'$ if $\cS(\cA)^f=\cS(\cA')$. The set of Cayley isomorphisms from $\cA$ to itself is denoted by $\iso_{cay}(\cA)$. Up to notation, the following statement is known as the Schur theorem on multipliers (see, e.g., \cite[statement~(1) of Theorem~2.3]{EP4}).

\lmml{070917a}
Let $\cA$ be an S-ring over an abelian group~$G$. Then 
$$
Z(\aut(G))\le\iso_{cay}(\cA).
$$
\elmm

The elements of $\cS$ and the number $\rk(\cA)=|\cS|$ are called, respectively, the {\it basic sets} and {\it rank} of the S-ring~$\cA$. Any union of basic sets is called an {\it $\cA$-subset of~$G$} or {\it $\cA$-set}; the set of all of them is denoted by $\cS(\cA)^\cup$. The
latter set is closed with respect to taking inverse and product. Given $X\in\cS(\cA)^\cup$ the submodule of~$\cA$ spanned by the set
$$
\cS(\cA)_X=\{Y\in\cS(\cA):\ Y\subset X\}
$$
is denoted by $\cA_X$.\medskip  

Any subgroup of $G$ that is an $\cA$-set is called an {\it $\cA$-subgroup} of~$G$ or {\it $\cA$-group}; the set of all of them is denoted by $\fG(\cA)$. With each $\cA$-set $X$, one can  naturally associate two $\cA$-groups, namely $\grp{X}$ and  $\rad(X)$. The following useful lemma was proved in \cite[Lemma 2.1]{EKP2}.

\lmml{260917i}
Let $\cA$ be an S-ring over a group $G$, $H\in\fG(\cA)$ and $X\in\cS(\cA)$. Then the cardinality of the set $X\cap Hx$ does not depend on $x\in X$.
\elmm

Let $S=U/L$ be a  section  of $G$. It is  called an  {\it $\cA$-section} if both  $U$ and $L$ are $\cA$-groups. Given $X\in\cS(\cA)_U$, the module
$$
\cA_S=\Span \{\pi_S(X):\ X\in\cS(\cA)_U\}
$$
is an S-ring over the group~$S$, where $\pi_S:U\to S$ is the natural epimorphism. The basic sets of $\cA_S$ are exactly the sets from the right-hand side of the formula.

\subsection{Wreath and tensor products}
Let $S=U/L$ be an $\cA$-section. The S-ring~$\cA$ is called an {\it $S$-wreath product} if $L\trianglelefteq G$ and $L\le\rad(X)$ for all basic sets $X$ outside~$U$; in this case, we write
\qtnl{050813a}
\cA=\cA_U\wr_S\cA_{G/L},
\eqtn
and omit $S$ when $U=L$. When the explicit indication of the section~$S$ is not important, we use the term {\it generalized wreath product} and omit~$S$ in the previous notation. The $S$-wreath product is {\it nontrivial} or {\it proper}  if $1\ne L$ and $U\ne G$.\medskip

If $\cA_1$ and $\cA_2$ are S-rings over groups $G_1$ and $G_2$ respectively, then the subring $\cA=\cA_1\otimes \cA_2$ of the ring $\mZ G_1\otimes\mZ G_2=\mZ G$, where $G=G_1\times G_2$, is an S-ring over the group $G$ with
$$
\cS(\cA)=\{X_1\times X_2: X_1\in\cS(\cA_1),\ X_2\in\cS(\cA_2)\}.
$$
It is called the {\it tensor product} of $\cA_1$ and $\cA_2$. The following statement was proved in \cite[Lemma 2.3]{EKP2}

\lmml{aux0}
Let $\cA$ be an S-ring over an abelian group $G=C\times D$. Assume that $C$ and $D$ are $\cA$-groups. Then 
\nmrt	
\tm{1}	$X_C$ and $X_D$ are basic sets of~$\cA$ for all $X\in \cS(\cA)$,
\tm{2}  $\cA\ge \cA_C\otimes \cA_D$, and the equality holds if $\cA_C$ or $\cA_D$ is the group ring.
\enmrt
\elmm

\subsection{Cyclotomic S-rings.}

An S-ring $\cA$ over a group $G$ is said to be {\it cyclotomic} if there exists $M\le\aut(G)$ such that
$$
\cS(\cA)=\orb(M,G).
$$
In this case, $\cA$ is denoted by $\cyc(M,G)$. Obviously, the group $K=G_{right}M$ satisfies condition~\eqref{031017u}. Thus, any cyclotomic S-ring is schurian. When $M$ is a multiplicative subgroup of a finite field~$\mF$, we say that $\cA$ is a cyclotomic S-ring over~$\mF$. For such a ring, the group $\iso_{cay}(\cA)$ contains the Frobenius automorphism of the field~$\mF$. The following statement is a special case of \cite[Theorem~5.1]{EP4})

\lmml{020917a}
For every prime $p$, each S-ring over a group $C_p$ is cyclotomic.
\elmm

\subsection{Duality.}
Let $G$ be an abelian group. Denote by $\wh G$ the group dual to~$G$, i.e., the group of all irreducible complex characters of~$G$. It is well known that there is a uniquely determined lattice antiisomorphism between the subgroups of~$G$ and~$\wh{G}$~\cite{S94}. The image of the group $H$ with respect to this antiisomorphism is denoted by~$H^\bot$.\medskip

For any S-ring $\cA$ over the group~$G$, one can define the dual S-ring $\wh{\cA}$ over~$\wh G$ as follows: two irreducible characters of $G$ belong to the same basic set of $\wh{\cA}$ if they have the same value on each basic set of~$\cA$ (for the exact definition, we refer to~\cite{EP4,EP16}). One can prove that 
$$
\rk(\wh{\cA})=\rk(\cA)
$$ 
and the S-ring dual to $\wh\cA$ is equal to~$\cA$. The following statement collect some facts on the dual S-rings proved in~\cite[Sec.~2.3]{EP16}.

\lmml{dual}
Let $\cA$ be an $S$-ring over an abelian group $G$. Then
\nmrt
\tm{1} the mapping $\fG(\cA)\to\fG(\wh{\cA})$, $H\mapsto H^\bot$ is a lattice antiisomorphism,
\tm{2}  $\wh{\cA_H}=\wh{\cA}_{\wh{G}/H^\bot}$ and $\wh{\cA_{G/H}}=\wh{\cA}_{H^\bot}$ for every $H\in \fG(\cA)$,
\tm{3} $\cA=\cyc(K,G)$ for $K\leq \aut(G)$ if and only if $\wh{\cA}=\cyc(K,\wh{G})$,
\tm{4} $\cA=\cA_1\otimes\cA_2$ if and only if $\wh{\cA}=\wh{\cA_1}\otimes\wh{\cA_2}$,
\tm{5} $\cA$ is the $U/L$-wreath product if and only if $\wh{\cA}$ is the $L^{\bot}/U^{\bot}$-wreath product.
\enmrt
\elmm

\subsection{Subdirect product.}

Let $U$ and $V$ be groups. Assume that $\varphi$ and $\psi$ be homomorphisms from $U$ and from $V$ onto isomorphic groups, i.e., there exists 
$$ 
f\in\iso(\im(\varphi),\im(\psi)).
$$
In this situation, one can define the subdirect product of the groups  $U$ and $V$ with respect to the homomorphisms $\varphi$, $\psi$, and the isomorphism $f$ as the following subgroup of $U\times V$:
$$
U\prod_{f}^{\varphi,\psi}V=
\{(u,v)\in U\times V:\ f(\varphi(u))=\psi(v)\}.
$$
It is easily seen that if the groups $\im(\varphi)$ and $\im(\psi)$ are trivial, then the subdirect product equals $U\times V$.\medskip 

In what follows, we are interested in the orbits of subdirect products of permutation groups, in which the homomorphisms $\varphi$ and $\psi$ are induced by the actions of the groups on imprimitivity systems associated with normal subgroups. More exactly, let $A$ and $P$ be groups, and we are given the following data:
\nmrt
\tm{P1} $U_0\trianglelefteq U\le\aut(A)$ and $V_0\trianglelefteq V\le\aut(P)$,
\tm{P2} $X_A\in\orb(U,A)$ and $X_P\in\orb(V,P)$,
\tm{P3} a bijection $f:\orb(U_0,X_A)\to\orb(V_0,X_P)$.
\enmrt
Under these conditions, the sets $\Pi_A=\orb(U_0,A)$ and $\Pi_P=\orb(V_0,P)$ form imprimitivity systems of the (transitive) groups 
$$
U^{X_A}\le\sym(X_A)\qaq V^{X_P}\le\sym(X_P).
$$
Denote by $\varphi$ and $\psi$ the natural epimorphisms from $U$ onto $U^{\Pi_A}$ and from $V$ onto~$V^{\Pi_P}$, respectively. Note that the permutation groups $U^{\Pi_A}$ and $V^{\Pi_P}$ are regular (this follows from the normality of the groups $U_0$ and $V_0$). Therefore, each isomorphism 
from $U^{\Pi_A}$ onto $V^{\Pi_P}$ is induced by a certain bijection of the form given in condition~(P3).

\lmml{030917u}
In the above notation, let $G=A\times P$. Assuming $f\in\iso(U^{\Pi_A},V^{\Pi_P})$, denote by $K$ the subdirect product of the groups $U$ and $V$ with respect to the homomorphisms $\varphi$, $\psi$, and $f$. Then $K\le\aut(A)\times\aut(P)$ and 
$$
\bigcup_{Y\in\Pi_A}Y\times Y^f\in\orb(K,G).
$$
\elmm

\section{S-rings over $E_9$}\label{020917f}

Up to Cayley isomorphism, there are exactly ten S-rings over $E_9$. This can be checked in a straightforward way or with the help of the GAP package COCO2~\cite{GAP}. In this section, we cite relevant properties of these S-ring to be used in Sections~\ref{130917a} and~\ref{270917a}. The first statement can be established by inspecting the above ten S-rings one after the other. 

\thrml{310817a}
Every S-ring over a group $E_9$  is Cayley isomorphic to one of the S-rings listed below:
\nmrt
\tm{1} $\cyc(M,\mF)$, where $\mF=\GF_9$ and $1<M\le\mF^\times$,
\tm{2} the tensor product of two S-rings over $C_3$,
\tm{3} the wreath product of two S-rings over $C_3$.
\enmrt
\ethrm

In statement~(1), $M$ is a cyclic group of order $2$, $4$, or $8$ (in the last case, the corresponding S-ring is of rank~$2$). In statement~(2), there are three S-rings of ranks~$4$, $6$, and $9$ (the last one is $\mZ E_9$). In statement~(3), there are four S-rings of ranks~$3$, $4$, $4$, and~$5$. To establish some properties of this S-rings, we need an auxiliary notion. \medskip

Let $G$ be an abelian group~$G$ and $X$ an orbit of a subgroup of $\aut(G)$; in particular, $X=X^{-1}$ or $X\cap X^{-1}=\varnothing$. A uniform partition\footnote{A partition of a set is said to be {\it uniform} if all its classes have the same cardinality} $\Pi$ of $X$ is said to be {\it regular} if the condition $X=X^{-1}$ implies that 
\nmrt
\tm{R1} $\Pi^{-1}=\Pi$,
\tm{R2} the permutation $Y\mapsto Y^{-1}$, $Y\in\Pi$, is either trivial or fixed point free,
\tm{R3} $\und{X}\circ\sum_{Y\in\Pi}\und{Y}\,\und{Y}^{-1}=\alpha\,\und{X}$ for some integer~$\alpha\ge 0$.
\enmrt
Thus, in the case $X\cap X^{-1}=\varnothing$, any uniform partition of $X$ is regular. If $X=X^{-1}$, then the partition of $X$ into one class is regular. A less trivial example is given by groups $L\trianglelefteq K\le\aut(G)$: in this case, one can take any $X\in\orb(K,G)$ and $\Pi=\orb(L,X)$.

\lmml{010917a}
Let $\cA$ be an S-ring over a group $G=E_9$, and let $X\in\cS(\cA)$. Then
for any regular partition $\Pi$ of $X$, there exists groups $L\le M\le\aut(G)$ such that
$$
\Pi=\orb(X,L)\qaq X\in\orb(M,G).
$$
Moreover, the group $M$ is cyclic unless $\cA$ is the tensor product of two  $S$-rings of rank~$2$ and $|X|=4$. In the latter case, $M=E_4$.
\elmm
\proof According to Theorem~\ref{310817a}, we have the following cases:
\nmrt
\tm{X1} $X$ is an orbit of a Singer subgroup $M$ of the group $\aut(G)\cong\GL(2,p)$; in particular, $|X|=2$, $4$, or~$8$;
\tm{X2} $X$ is an orbit of a subgroup of  $\aut(C)\times\aut(C')\le\aut(G)$, where $C\cong C_3\cong C'$ are such that $G=C\times C'$; in particular, $|X|=1$, $2$, or~$4$;
\tm{X3} $X$ is an orbit of a subgroup $M\le\aut(G)$ of order~$6$ that stabilizes a group $C\cong C_3$; in particular, $|X|=1$, $2$, $3$, or~$6$.
\enmrt
Let $\Pi$ be a regular partition of~$X$. Without loss of generality, we may assume that $1<|\Pi|<|X|$. Since also $|\Pi|$ divides $|X|$, 
\qtnl{201017a}
(|X|,|\Pi|)=(4,2),\ (8,2),\ (8,4),\ (6,2),\ \text{or}\ \,(6,3),
\eqtn
where the first pair appears in cases (X1) and (X2), the second two appear in (X1), and the last two appear in~(X3). In all these cases, the set $X$ is symmetric. A simple counting argument using conditions (R2) and (R3) shows that the permutation defined in (R2) must be trivial unless $(|X|,|\Pi|)=(4,2)$ for the case (X2), and $(|X|,|\Pi|)=(6,2)$. Therefore, the number of possible partitions $\Pi$ of the set~$X$ is $1$ or~$3$, $3$, $1$, $1$, and~$1$, respectively to cases listed in~\eqref{201017a}. Now a straightforward check in each case completes the proof.\eprf\medskip

Lemma~\ref{010917a} shows that the group $L$ equals the kernel of the  homomorphism from~$M$ to $\sym(\Pi)$ induced by the action of $M$ on the $\orb(X,L)$. In what follows, we say that $(M,L)$ is a {\it standard pair} for the basic set~$X$ and regular partition~$\Pi$; though the standard pair  is not uniquely determined, the following statement holds true.

\lmml{030917a}
In the notation of Lemma~\ref{010917a}, assume that the group $M^\Pi$ is cyclic. Then for any regular cyclic group $C\le\sym(\Pi)$ centralizing the permutation in condition~{\rm (R2)}, there exist $\sigma\in\aut(G)$ such that 
\nmrt
\tm{1} $(M^\sigma,L)$ is a standard pair for $X$ and $\Pi$, 
\tm{2} $(M^\sigma)^\Pi=C$.
\enmrt
\elmm
\proof The statement is trivial if $|\Pi|\le 3$, because $\sym(\Pi)$ contains a unique cyclic subgroup of order~$|\Pi|$. Since $|X|\le 8$ and $|\Pi|$ divides $|X|$, we may assume that
$$
(|X|,|\Pi|)=(4,4),\ (6,6),\ (8,8),\ \text{or}\ \,(8,4).
$$
The condition on $C$ implies that $C\le\aut(G)$ in the first three cases. This proves the required statement in these cases, because any two cyclic subgroups of $\aut(G)\cong\GL(2,3)$ of the same order at least $4$ are conjugate. In case $(8,2)$, the condition on $C$ implies that permutation in condition~{\rm (R2)} belongs to~$C$.  This leaves exactly three possibilities for $C$, and for each of them $C=M^\Pi$, where $M$ is one of the three Singer subgroups of~$\aut(G)$.\eprf

\section{Dense S-rings over $E_9\times C_p$: basic sets}\label{130917a}

Throughout this section, we assume that $p>3$ and $G=A\times P$, where $A=E_9$ and $P=C_p$. In what follows, $\cA$ is a {\it dense} $S$-ring over~$G$, which means that $A$ and~$P$ are $\cA$-subgroups of~$G$. By Lemma~\ref{070917a}, we have
\qtnl{290917y}
\grp{\tau}\times\aut(P)=Z(\aut(G))\le\iso_{cay}(\cA),
\eqtn 
where $\tau\in\aut(A)$ is the involution taking $a$ to $a^{-1}$.\medskip 


Let $X$ be a basic set of the S-ring~$\cA$. In view of Lemma~\ref{aux0}, the projections $X_A\subset A$ and $X_P\subset P$ are basic sets of the S-rings $\cA_A$ and $\cA_P$. Therefore 
$$
X_A\times X_P\in\cS(\cA)^\cup. 
$$
By Lemmas~\ref{010917a} and~\ref{020917a}, there exist groups $U(X)\le\aut(A)$ and $V(X)\le\aut(P)$ such that
$$
X_A\in\orb(U,A)\qaq X_P\in\orb(V,P),
$$
where $U=U(X)$ and $V=V(X)$. For any element $a\in X_a$,  each basic set inside $X_A\times X_P$ intersects $\{a\}\times X_P$. Since the group $V\le\aut(G)$ acts transitively on the latter set, formula~\eqref{290917y} implies that $V$ acts transitively on $\cS(\cA)_{X_A\times X_P}$.

\lmml{020917d}
Let $\Pi_P(X)=\{X(a)\}_{a\in X_A}$, where
\qtnl{130917c}
X(a)=\{x\in X_P:\ (a,x)\in X\}.
\eqtn
Then there exists a group $V_0\le V$ such that
$$
\Pi_P=\orb(V_0,X_P).
$$
In particular, $\Pi_P$ is an imprimitivity system for the group $V^{X_P}$.
\elmm
\proof Denote by $V_0$ the subgroup of $V$ leaving $X$ fixed (as a set). From formula~\eqref{290917y}, it follows that each set $X_a$ is contained in some $Y\in\orb(V_0,X_P)$. On the other hand, $X(a)$ cannot be smaller than $Y$ by the definition of $V_0$. Thus, $X_a=Y$.\eprf\medskip

Let us define an equivalence relation $\sim$ on the set $X_A$ by setting $a\sim b$ if and only if $X(a)=X(b)$. In particular, all the elements of $X_A$ are $\sim$-equivalent if and only if $X=X_A\times X_P$. Denote by $\Pi_A=\Pi_A(X)$ the partition of $X_A$ into the classes of the equivalence relation~$\sim$. From our definitions and Lemma~\ref{020917d}, it follows that the mapping
\qtnl{030917av}
f:\Pi_A\to\Pi_P,\ [a]\mapsto X(a),
\eqtn
is a well-defined bijection, where $[a]$ denotes the class of the equivalence relation~$\sim$ that contains $a\in X_A$. Moreover,
\qtnl{030917d}
X=\bigcup_{Y\in \Pi_A}Y\times Y^f
\eqtn

\lmml{020917e}
$\Pi_A$  is a regular partition of $X$ in the sense of Section~\ref{020917f}. 
\elmm
\proof It is easily seen that $[a]y=X\cap Ay$ for all $a\in X_A$ and $y\in X(a)$. By Lemma~\ref{260917i}, this implies that the partition $\Pi_A$ is uniform. Without loss of generality, we may assume that $X_A$ is symmetric. By formula~\eqref{290917y} and since $\aut(P)$ acts transitively on  $\cA_{X_A\times X_P}$, there exists $\sigma\in V$ such that $(X^\tau)^\sigma=X$. In view of equality~\eqref{030917d} and Lemma~\ref{020917d}, this implies that 
$$
([a]\times X(a))^{\tau\sigma}=[a]^{-1}\times X(a)^\sigma=[b]\times X(b)
$$
for all $a\in X_A$, where the element $b\in X_P$ is defined by the condition~$X(a)^\sigma=X(b)$. Thus, $[a]^{-1}=[b]$ and the partition $\Pi_A$ satisfies condition~(R1). Furthermore, if $[a]=[a]^{-1}$ for some $a\in X_A$, then $X^\tau=X$ by formula~\eqref{030917d}, and hence the permutation $[a]\to[a]^{-1}$ is trivial. This shows that $\Pi_A$ satisfies condition~(R2). Finally, again by formula~\eqref{030917d} we have
$$
\sum_{Y\in\Pi_A}\und{Y}\,\und{Y}^{-1}=\und{A}\circ(\und{X^{}}\,\und{X}^\tau)=\alpha\und{X_A}+\xi
$$
for some integer $\alpha\ge 0$ and $\xi\in\cA$ such that $\und{X_A}\circ\xi=0$. This shows that $\Pi_A$ satisfies condition~(R3).\eprf\medskip

From Lemmas~\ref{020917e} and~\ref{010917a}, it follows that given $X\in\cS(\cA)$ there exists a standard pair~$(U,U_0)$ for the set~$X_A\in\cS(\cA)$ and regular partition~$\Pi_A$. The following statement is the main result of this section.

\thrml{set}
In the above notation, the set $X$ is an orbit of the subdirect product $K=K(X)$ of the groups $U$ and $V$ with respect to the homomorphisms 
$$
\varphi:U\to U^{\Pi_A}\qaq\psi:V\to V^{\Pi_P},
$$ 
and the isomorphism~$f:U^{\Pi_A}\to V^{\Pi_P}$ induced by bijection~\eqref{030917av}. Furthermore, $K\le\aut(G)$
\ethrm
\proof According to our notation, we are in the situation described by the conditions (P1), (P2), and~(P3). Moreover, the group $V^{\Pi_P}$ is cyclic. First, assume that the group $U^{\Pi_A}$ is also cyclic. Then by Lemma~\ref{030917a} for $(M,L)=(U,U_0)$, $\Pi=\Pi_A$, and $C=f^{-1}V^{\Pi_P}_Pf$, the standard pair can be chosen so that 
$$
U^{\Pi_A}=f\,U^{\Pi_P}f^{-1},
$$
i.e., $f\in\iso(U^{\Pi_A},V^{\Pi_P})$. Thus, from Lemma~\ref{030917u} and relation~\eqref{030917d} it immediately follows that
$$
K\le\aut(A)\times\aut(P)=\aut(G)\qaq X\in\orb(K,G),
$$
as required.\medskip

To complete the proof, we show that the group $U^{\Pi_A}$ must be cyclic. Assume on the contrary that this is not true. Then $|\Pi_A|= 4$. Moreover, by statement~(1) of Lemma~\ref{010917a}, the S-ring~$\cA_A$ is the tensor product of two trivial $S$-rings and $|X|=4$. In particular, there are two distinct $\cA_A$-groups $C$ and $D$, each of order~$3$ and
$$
X=C^\#\times D^\#.
$$
Note that both $C$ and $D$ are also $\cA$-groups, and  $G=C\times(DP)=D\times(CP)$. By statement~(1) of Lemma~\ref{aux0}, this implies that $Y=X_{DP}$ and $Z=X_{CP}$ are basic sets of~$\cA$. It is easily seen that 
$$
Y_P=Z_P=X_P\qaq |\Pi_P(Y)|=|\Pi_P(Z)|=2\qaq \Pi_P(Y)\ne \Pi_P(Z).
$$
By Lemma~\ref{020917d}, this implies that the transitive cyclic group $V(X)^{X_P}$ has two distinct imprimitivity systems, each with exactly two blocks, a contradiction.\eprf\medskip

For distinct basic sets $X$ and $Y$ of the S-ring~$\cA$, the groups $K(X)$ and $K(Y)$ are not necessarily equal: even if the standard pairs for $X_A$ and $Y_A$ are equal,  the subdirect products $K(X)$ and $K(Y)$ may correspond to different isomorphisms~$f$. The following statement provide a sufficient condition for $Y$ to be an orbit of $K(X)$. In what follows, we set
$$
\fG(\cA)'=\{H\in\fG(\cA):\ G=H\times H'\ \,\text{for some}\ \, H'\in\fG(\cA)\}.
$$

\lmml{100917u}
Let $\cA$ be a dense $S$-ring over~$G$, $X,Y\in\cS(\cA)$, and $K=K(X)$. Then $Y\in\orb(K,G)$ if at least one of the following conditions is satisfied:
\nmrt
\tm{1} $Y=X_H$ for some $H\in\fG(\cA)'$,
\tm{2} $Y=X^\sigma$ for some $\sigma\in\ Z(\aut(G))$,
\tm{3} $X_A=Y_A$. 
\enmrt
\elmm
\proof Under condition~(1), the $\cA$-groups $H$ and $H'$ are $K$-invariant. Therefore, it is easily seen that
$$
X\in\orb(K,G)\quad\Rightarrow\quad X_{H^{}},X_{H'}\in\orb(K,G),
$$
and we are done. Now assume that condition~(2) is satisfied. Since the automorphism $\sigma$ centralizes $K\le\aut(G)$, we conclude that 
$$
Y=X^\sigma\in\orb(K^\sigma,G)=\orb(K,G),
$$
as required. To complete the proof, it suffices to note that condition~(3) is a consequence of conditions~(1) and~(2).\eprf

\section{Dense S-rings over $E_9\times C_p$ are cyclotomic}\label{270917a}

The main result of the present section is given in the following theorem. Along the proof, we freely use the notation introduced in Section~\ref{130917a}

\thrml{240917a}
Every dense S-ring over~$E_9\times C_p$ is cyclotomic.
\ethrm
\proof Let $\cA$ be a dense S-ring over the group~$G=A\times P$ with $A=E_9$ and $P=C_p$ for a prime~$p>3$ (for $p=2,3$, the required statement can be verified by enumeration of the S-rings over small groups~\cite{Ziv}). We divide the proof into three separate cases depending on which statement of Theorem~\ref{310817a} holds for the S-ring~$\cA_A$.\medskip

{\bf Case 1:} $\cA_A=\cyc(M,\mF)$, where $\mF=\GF_9$ and $1<M\le\mF^\times$.
In this case, $M$ is a cyclic group of order $m\in\{2,4,8\}$. Fix a basic set $$
X\in\cS(\cA)_{G\setminus(A\cup P)},\ |X_A|=m.
$$
Denote by $K$ the group $K(X)$ defined in Theorem~\ref{set}. First assume that $|M|\ne 4$. Then $\rk(\cA_A)=2$ or $\fG(\cA_A)$ contains all subgroups of~$A$. It easily follows that 
\qtnl{201017a1}
\cS(\cA)^\#=\{(X_H)^\sigma:\ H\in\fG'(\cA),\ \sigma\in Z(\aut(G))\}.
\eqtn
By  Lemma~\ref{100917u}, this implies that $\cA=\cyc(K,G)$.\medskip

Now let $|M|=4$. In this case, for any $Y\in\cS(\cA)^\#$, the set $Y_A$ is either trivial, or is equal to~$X_A$ or to~$A^\#\setminus X_A$. By Lemma~\ref{100917u}, it suffices to verify that in the last case, $X$ and $Y$ are orbits of a certain group $K'\le\aut(G)$ (except for one case, $K'$ will be equal to~$K$). To this end, set
\qtnl{210917a}
k_X=|\Pi_A(X)|\qaq k_Y=|\Pi_A(Y)|.
\eqtn
Each of the numbers $k_X$ and $k_Y$ divides $|X|=|Y|=4$, and hence is equal to~$1$, $2$, or~$4$.  Let us analyze all these possibilities. It is convenient to denote the eight nontrivial elements of the group~$A$ by $a_i^{\pm 1}$, $i=1,2,3,4$, so that the orbits $X_A$ and $Y_A$ of the group~$M$ are of the form:
$$
X_A=\{a^{}_1,a^{}_3,a_1^{-1},a_3^{-1}\}\qaq
Y_A=\{a^{}_2,a^{}_4,a_2^{-1},a_4^{-1}\}.
$$
Without loss of generality, we may assume that $a^{}_2=a^{}_1a^{}_3$ and $a^{}_4=a^{}_1a_3^{-1}$, and also 
\qtnl{021017a}
f_X([a_1])=f_Y([a_2])\quad\text{and hence}\ \, f_X([a^{}_1])=f_Y([a^{}_2]),
\eqtn
where $f_X$ and $f_Y$ are the bijections defined by formula~\eqref{030917av} for $X$ and $Y$, respectively.\medskip

{\bf Claim~1:} $\{k_X,k_Y\}\ne\{4,2\}$ and  $\{k_X,k_Y\}\ne\{4,1\}$. Assume, for instance, that $k_X=4$. Then a straightforward calculation shows that
\qtnl{130917d}
\begin{split}
(X_A\,X)\cap(Y_A\times X_P)=&
a^{\phantom{-1}}_2X(a^{\phantom{-1}}_1,a^{\phantom{-1}}_3)\,\cup\,
a^{\phantom{-1}}_4X(a^{\phantom{-1}}_1,a^{-1}_3)\,\cup\,\\
& a^{-1}_2X(a^{-1}_1,a^{-1}_3)\,\cup\,
a^{-1}_4X(a^{-1}_1,a^{\phantom{-1}}_3),
\end{split}
\eqtn
where $X(a_i,a_j)=X(a_i)\cup X(a_j)$ for all~$i,j$. Note that the left-hand side of~\eqref{130917d} is an $\cA$-set, because $X_A$, $X$, $Y_A$, and $X_P$ are basic sets of~$\cA$. Furthermore, assumption~\eqref{021017a} implies that it contains $a_2X(a_1)$ and hence intersects~$Y$ nontrivially. Thus,
$$
Y\subseteq (X_A\,X)\cap(Y_A\times X_P).
$$
On the other hand, from the form of the right-hand side of~\eqref{130917d}, it follows that $k_Y\ne 1$, and if $k_Y=2$, then the cyclic group $V(Y)=V(X)$ has two different imprimitivity systems, each consisting two blocks (Lemma~\ref{020917d}). Since this is impossible, the claim is proved.\eprf\medskip

{\bf Claim~2:} if $\{k_X,k_Y\}=\{1,2\}$, then $X,Y\in\orb(K',G)$ for some group $K'$ such that $K<K'\le\aut(G)$.  Without loss of generality, we may assume that $k_X=1$ and $k_Y=2$. Note that the Frobenius automorphism $k'$  of the field~$\mF$ is an automorphism of the S-ring~$\cA_A$. It follows that $X_A,Y_A\in\orb(U',A)$, where $U'=\grp{U,k'}$. Set 
\qtnl{220917k}
K'=U'\prod_{f}^{\varphi',\psi}V,
\eqtn
where $\varphi':U'\to U'/U'_0$ and $U'_0=\grp{M_0,k'}$ with $M_0$ being the subgroup of $M$ of order~$2$. Then by Lemma~\ref{030917u} applied to the set~$X$ and trivial bijection $f:\{X_A\}\to\{X_P\}$, and to the set~$Y$ with the bijection $f_Y$, we conclude that  $X$ and $Y$ are orbits of the group~$K'\le\aut(G)$.\eprf\medskip

By Claims 1 and 2, to complete the proof of the Case~1, we may assume that $k_X=k_Y:=k$. If now $k=1$ or $2$, then the bijection~\eqref{030917av} is unique and hence $Y\in\orb(K,G)$. Assume that $k=4$. In this case, the groups $K(X)$ and $K(Y)$ may correspond to subdirect products with different bijections $f_X$ and $f_Y$. Namely, there are two possibilities:
$$
f_X([a^{}_3])=f_Y([a^{}_4])\qoq f_X([a^{}_3])=f_Y([a_4^{-1}]).
$$
However, the first case is impossible, because the $\cA$-set  
$$
(XY^{-1})\cap A^\#=\{a_2^{\pm 1},a_3^{\pm 1}\}
$$ 
intersects each of the two different basic sets $X_A$ and $Y_A$ nontrivially.
Since in the last case, $Y\in\orb(K,G)$, we are done.\medskip

{\bf Case 2:} $\cA_A=\cA_C\otimes\cA_D$, where $C$ and $D$ are subgroups of $A$ such that $A=C\times D$ and $|C|=|D|=3$. First assume that one of the S-rings $\cA_C$ or $\cA_D$ is the group ring, say the first one. Then by statement~(2) of Lemma~\ref{aux0} for $G_1=C$ and $G_2=DP$, we have
\qtnl{100917uv}
\cA=\mZ C\otimes \cA_{DP}.
\eqtn
Since $DP$ is a cyclic group and $D,P$ are $\cA_{DP}$-groups, the classification of S-rings over a cyclic group $C_{pq}$ with primes~$p$ and~$q$  implies that the S-ring  $\cA_{DP}$ is cyclotomic (see~ \cite{KP81}). Since also the S-ring $\mZ C$ is cyclotomic, formula~\eqref{100917uv} shows that so is the S-ring~$\cA$.\medskip

To complete this case assume that both $\cA_C$ and $\cA_D$ are S-rings of rank~$2$. Take any $X\in\cS(\cA)_{G\setminus A}$ such that $X_A=C^\#\times D^\#$. Then one can easily verify that formula~\eqref{201017a1} holds. By Lemma~\ref{100917u}, this implies that $\cA=\cyc(K,G)$ with $K=K(X)$.\medskip

{\bf Case 3:} $\cA_A=\cA_C\wr\cA_{A/C}$, where $C\le A$ is a group of order~$3$. Depending on whether $\cA_{A/C}$ is of rank~$2$ or $3$, the set $\cS(\cA)_{A\setminus C}$ consists of one set of cardinality~$6$  or two sets of cardinalities~$3$, respectively. Fix a basic set
$$
X\in\cS(\cA)_{(A\setminus C)\times P^\#}.
$$
First, assume that $\cA_C$ is of rank~$3$. Since the number $|\Pi_A(X)|$ divides~$6$, there exists a standard pair $(U,U_0)$ for $X_A$ and $\Pi_A(X)$ such that the set $\orb(U_0,C)$ consists of singletons. Assume that the group~$K=K(X)$ is associated with this pair. Then obviously
$$
Y\in\orb(K,G)\quad\text{for all}\ \,Y\in\cS(\cA)_{C\times P}.
$$
Next, we observe that $X_A\times X_P$ is the union of $X$ and $X^{-1}$,  
where depending on whether $\cA_{A/C}$ is of rank~$2$ or~$3$, we have
$$
X=X^{-1}\ \,\text{and}\ \,|X_A|=6\qoq 
X\cap X^{-1}=\varnothing\ \,\text{and}\ \,|X^{}_A|=|X^{-1}_A|=3,
$$
respectively. This easily implies that any $Y\in \cS(\cA)$ contained in $(A\setminus C)\times P$ is of the form $X^\sigma$ for some $\sigma\in Z(\aut(G))$. Thus, again $Y\in \orb(K,G)$ by Lemma~~\ref{100917u} and hence $\cA=\cyc(K,G)$.\medskip

Let now $\cA_C$ and $\cA_{A/C}$ be of rank~$2$ and $3$, respectively. Then
\qtnl{220917d}
A^\bot,P^\bot\in\fG(\wh{\cA})
\eqtn
by statement~(1) of Lemma~\ref{dual}. It follows that $\wh{\cA}$ is a dense S-ring over the group~$\wh G$. The statements~(2) and~(3) of that lemma imply that the restriction of $\wh \cA$ to the group~$A^\bot$ is the wreath product of the S-rings of rank~$3$ and rank~$2$. By the previous paragraph, we conclude that $\wh{\cA}$ is a cyclotomic S-ring over~$\wh{G}$. By statement~(3) of Lemma~\ref{dual}, this proves that $\cA$ is a cyclotomic S-ring over~$G$.\medskip

In the remaining case, both $\cA_C$ and $\cA_{A/C}$ are of rank~$2$. It follows that $\cA_A$ is of rank~$3$ and 
$$
\cS^\#(\cA)=\{C^\#,A\setminus C\}.
$$
Fix arbitrary basic sets $X,Y\in\cS(\cA)$ such that 
$$
X_A=A\setminus C,\quad Y_A=C^\#,\quad X_P=Y_P\ne 1_P.
$$
By Lemma~\ref{100917u}, it suffices to verify that  $X$ and $Y$ are orbits of a certain group $K'\le\aut(G)$. To this end, we define the numbers $k_X$ and $k_Y$ by formula~\eqref{210917a}. Then
from Lemmas~\ref{020917d} and~\ref{020917e}, it follows that
\qtnl{220917j}
k_X\in\{1,2,3,6\}\qaq k_Y\in\{1,2\}.
\eqtn
As in Case~1, not each combination for the pair $(k_X,k_Y)$ is possible.\medskip

{\bf Claim 3:} $(k_X,k_Y)\ne (6,1)$ and $(k_X,k_Y)\ne (3,2)$. Let us consider the first case. Denote by $r$ the cardinality of $X_P$. Then 
$$
|X|=r\qaq |Y|=2r.
$$
The set $G\setminus (A\cup P)$ is partitioned into basic sets $X^\sigma$ and $Y^\sigma$, where $\sigma\in Z(\aut(G))$. Since $|X^\sigma|=|X|$, $|Y^\sigma|=|Y|$, and $|X_P|=|Y_P|=(p-1)/r$, we obtain
\qtnl{220917a}
|\cS(\cA)|=|\cS(\cA_A)|+|\cS(\cA_P)|-1+\frac{7(p-1)}{r}.
\eqtn
Next, let $\pi:G\to G/C$ be the natural epimorphism. Since $k_X=6$, each of the sets $\pi(X^\sigma)$ is of cardinality~$r$. It follows that 
\qtnl{220917c}
|\cS(\cA_{G/C})|=|\cS(\cA_{A/C})|+|\cS(\cA_{CP/C})|-1+\frac{2(p-1)}{r}.
\eqtn
For the S-ring $\wh{\cA}$ dual to $\cA$, relation~\eqref{220917d} holds. Since the S-rings $\cA_A$ and $\cA_{G/A}$ as well as $\cA_P$ and $\cA_{G/P}$ are isomorphic, statement~(2) of Lemma~\ref{dual} and  equality~\eqref{220917a} yield
\qtnl{220917f}
|\cS(\wh{\cA})_{\wh G\setminus (A^\bot\cup P^\bot)}|=\frac{7(p-1)}{r}.
\eqtn
Furthermore, $C^\bot\cong C_{3p}$ is an $\wh{\cA}$-group and the restriction of $\wh{\cA}$ to this group is isomorphic to the S-ring $\cA_{G/C}$. Therefore from equality~\eqref{220917c}, it follows that
\qtnl{220917g}
|\cS(\wh{\cA})_{C^\bot\setminus(A^\bot\cup P^\bot)}|=\frac{2(p-1)}{r}.
\eqtn
Now using equalities~\eqref{220917f} and \eqref{220917g}, we conclude that
there are exactly $5(p-1)/r$ basic sets of $\wh{\cA}$ outside $A^\bot$, $P^\bot$, and $C^\bot$. All of these basic sets are obtained from any one of them by applying an automorphism from $Z(\aut(\wh G))$. Therefore they have the same size, say~$k$. This implies that 
\qtnl{220917i}
k\frac{5(p-1)}{r}=|\wh{G}\setminus(A^\bot\cup P^\bot\cup C^\bot)|=6(p-1).
\eqtn
On the other hand, the S-rings $\cA_P=\cyc(V,P)$ and $\wh{\cA}_{P^\bot}$ are isomorphic, where $V=V(X)=V(Y)$ is a subgroup of $\aut(P)$ of order~$r$. Therefore the nonidentity basic sets of the latter S-ring are of cardinality~$r$. It follows that $r$ divides~$k$, which contradicts equality~\eqref{220917i}.\medskip

The proof of the Claim~3 for the case  $(k_X,k_Y)=(3,2)$ differs from the previous argument only in the values of the parameters. Namely, here
$|X|=2r$ and $|Y|=r$. Therefore, the last summands on the right-hand sides in formulas~\eqref{220917a} and~\eqref{220917f} are equal to $5(p-1)/r$, and those in formulas~\eqref{220917c} and~\eqref{220917g} are $(p-1)/r$. Thus, equality \eqref{220917i} leads to the equality $4k/r=6$, which is also impossible, because $r$ divides~$k$. The claim is proved.\medskip\eprf

Let us return to the remaining part of Case~$3$. By Claim~3 and formula~\eqref{220917j}, there are six possibilities for the pair $(k_X,k_Y)$. For each of them, we define a group $K'\le\aut(G)$ by formula~\eqref{220917k}, where the standard pair $(U',U'_0)$ is given in the second and third columns of Table~\ref{tbl1} below; in the fourth column contains the sizes of the $U'_0$-orbits.

\begin{table}[h]
\begin{center}
\begin{tabular}{|c|c|c|c|}
\hline
$(k_X,k_Y)$ & $U'$      &  $U'_0$    & $\orb(U'_0,A)$\\
\hline
$(1,1)$     &  $C_6$    &  $C_6$     & $[1,2,6]$  \\
$(1,2)$		&  $D_{12}$ &  $\sym(3)$ & $[1,1,1,6]$ \\
$(2,1)$     &  $D_{12}$ &  $\sym(3)$ & $[1,2,3,3]$ \\
$(2,2)$     &  $C_6$    &  $C_3$     & $[1,3,3,1,1]$\\
$(3,1)$	    &  $C_6$    &  $C_2$     & $[1,2,2,2,2]$\\
$(6,2)$	    &  $C_6$    &  $1$       & $[1,\ldots,1]$  \\
\hline
\end{tabular}
\medskip
\end{center}
\caption{Standard groups for Case~3}
\label{tbl1}
\end{table}

\noindent A straightforward check shows that in each case, $Y\in\orb(K',G)$, as required.\eprf

\section{The proof of Theorem~\ref{230917a}}\label{031017i}

We deduce Theorem~\ref{230917a} in the end of the section from the theorem below giving a complete description of all S-rings over the group $E_9\times C_p$.

\thrml{main}
Let $\cA$ be an $S$-ring over a group $G=E_9\times C_p$, where $p>3$ is a prime. Then one of the following statements holds:
\nmrt
\tm{1} $\cA$ is trivial or cyclotomic,	
\tm{2} $\cA$ is the tensor product of a trivial S-ring and an S-ring over $C_3$,
\tm{3} $\cA$ is a proper $S$-wreath product with $|S|\leq 3$.
\enmrt
\ethrm
\proof If the S-ring $\cA$ is dense, then we are done with statement~(1) by Theorem~\ref{240917a}. Assume that $\cA$ is not dense. Then $A$ or $P$ is not an $\cA$-subgroup of $G$. By the duality (see Lemma~\ref{dual}), we may assume that $A\not\in\fG(\cA)$. Denote by $C$ the maximal $\cA$-subgroup of the group~$A$. Clearly, this group is trivial or of order~$3$. The lemma  below is a special case of \cite[Lemma 6.2]{EKP2}.

\lmml{aux2}
In the above notation, one of the following statements holds:
\nmrt
\tm{1} $\cA=\cA_C\wr \cA_{G/C}$ and also $\rk(A_{G/C})=2$,
\tm{2} $\cA$ is the $U/L$-wreath product, where $P\leq L<G$ and $U=CL$.
\enmrt
\elmm

Without loss of generality, we may assume that $\cA$ is as in statement~(2) of Lemma~\ref{aux2}: otherwise this lemma implies that either $\cA$ is trivial (if $C=1$) and statement~(1)  of Theorem~\ref{main} holds, or  statement~(3) of this theorem holds with $S=C/C$. Furthermore, if $U<G$ then $|U/L|\leq 3$ and statement~(3) of Theorem~\ref{main} holds. Thus, we may also assume that $U=G$. Then $C\not\leq L$, for otherwise  $G=L$, a contradiction. Since $|C|=3$, it follows that $C\cap L=e$. Thus,
$$
G=C\times L\qaq |L|=3p.
$$
Note that the S-ring $\cA_L$ is circulant. Moreover, since $A$ is not an $\cA$-group, the subgroup of $L$ of order~$3$ is not an $\cA_L$-group. According to \cite{KP81}, this implies that
\qtnl{260917u}
\rk(\cA_L)=2\qoq \cA_L=\cA_P\wr\cA_{L/P}.
\eqtn
Assume that $\cA_C=\mZ C$. Then  $\cA=\cA_C\otimes\cA_L$ by statement~$2$ of Lemma~\ref{aux0}. In particular, statement~(2) of Theorem~\ref{main} holds, whenever $\rk(\cA_L)=2$. On the other hand, if $\cA_L$ is not trivial, then $\cA$ is obviously the $CP/P$-wreath product and statement~(3) of Theorem~\ref{main} holds. Thus, we may assume that
$$
\rk(\cA_C)=2.
$$
Denote by $L_0$ the trivial subgroup of~$L$ if $\rk(\cA_L)=2$, and the group $P$ otherwise. In view of~\eqref{260917u}, we have $L_0\in\fG(\cA)$. In particular, $L\setminus L_0$ is an $\cA$-set.

\lmml{260917v}
If $X\in\cS(\cA)$ is contained in $C^\#\times (L\setminus L_0)$, then $X=X_C\times X_L$.
\elmm
\proof For all $X\in\cS(\cA)$ contained in $C^\#\times L^\#$, we have 
$$
|X_L|\le |X|\le |C^\#|\,|X_L|=2|X_L|.
$$ 
It follows that $X_C\times X_L$ is the union of two basic sets $X$ and $X'$ of the same cardinality. Now if $X=X'$, then $X=X_C\times X_L$ and we are done. In the remaining case, 
\qtnl{260917o}
|X|=\frac{|X_C|\cdot|X_L|}{2}=\frac{2|X_L|}{2}=|X_L|. 
\eqtn
Assume that $X_L\subseteq L\setminus L_0$. From the definition of the group~$L_0$, it follows that $|L\setminus L_0|\le 3p-1$. Therefore, equality~\eqref{260917o} yields
\qtnl{260917p}
|X|=|X_L|\le|L\setminus L_0| \le 3p-1.
\eqtn
On the other hand, if $L_0=1$, then $\rk(\cA_L)=2$ and hence $X_L=L\setminus L_0$. Furthermore, if $L_0=P$, then $\cA_L=\cA_P\wr\cA_{L/P}$ and hence $X_L=P$ or $X=L\setminus P$.  However, the first case is impossible, because
by Lemma~\ref{260917i} for $H=P$ the number $|X|$ must be even. Thus, in any case, $X_L=L\setminus L_0$ and hence
$$
X^{}\,\cup\, X'=C^\#\times (L\setminus L_0).
$$
The the right-hand side includes the set $C_0:=A\setminus C$ of cardinality~$6$. Therefore, at least one of $X$ or $X'$, say $X$, contains three elements from $C_0$. According to \cite[Lemma 6.1]{EKP2}, this implies that 
$$
|X|\ge |(X\cap C_0)P|\ge 3p,
$$
which contradicts inequality~\eqref{260917p}.\eprf\medskip

From Lemma~\ref{260917v}, it follows that if $\rk(\cA_L)=2$, then $\cA=\cA_C\otimes\cA_L$ and we are done with statement~(2) of Theorem~\ref{main}. To complete the proof, in view of~\eqref{260917u} we may assume that $\cA_L=\cA_P\wr\cA_{L/P}$. In this case, statement~(1) of Lemma~\ref{aux0} and Lemma~\ref{260917v} imply that 
$$
P\le\rad(X)\quad\text{for all}\ \,X\in\cS(\cA)_{G\setminus CP}.
$$
It follows that $\cA$ is the $CP/P$-product and we and we are done with statement~(3) of Theorem~\ref{main}.\eprf\medskip

{\bf Proof of the Theorem~\ref{230917a}.} For $p\leq 3$, the statement  follows from the computational results obtained in~\cite[p.~498]{Ziv}. Let $p>3$ and $\cA$ an $S$-ring over~$G$. Then by Theorem~\ref{main}, this ring is obviously schurian if statement~(1) of this theorem holds. In case of statement~(2), the S-ring $\cA$ being the tensor product of two schurian S-rings is also schurian. To complete the proof, we may assume that $\cA$ is a proper $S$-wreath product with $|S|\leq 3$. Then the S-ring $\cA_S$ is either trivial or a group ring. Thus, the group $\aut(\cA_S)$ is permutation isomorphic to $S_{right}$, $\sym(3)$, or $\alt(3)$. In any case, according to a criterion of schurity of a generalized wreath product \cite[Corollary 10.3]{MP3}, the $S$-ring~$\cA$ is schurian.\eprf

\end{document}